


\magnification=\magstep 1
\parindent = 0 pt
\baselineskip = 16 pt
\parskip = \the\baselineskip

\font\AMSBoldBlackboard = msbm10

\def\RR{{\hbox{\AMSBoldBlackboard R}}}
\def\CC{{\hbox{\AMSBoldBlackboard C}}}
\def\AA{{\hbox{\AMSBoldBlackboard A}}}
\def\QQ{{\hbox{\AMSBoldBlackboard Q}}}
\def\PP{{\hbox{\AMSBoldBlackboard P}}}
\def\ZZ{{\hbox{\AMSBoldBlackboard Z}}}

\settabs 12\columns
\rightline{math.NT/9809119}
\vskip 1 true in
{\bf \centerline{THE EXPLICIT FORMULA AND A PROPAGATOR}}
\vskip 0.5 true in
\centerline{Jean-Fran\c{c}ois Burnol}
\par
\centerline{September 1998}
\centerline{revised November 1998}
\par
I give a new derivation of the Explicit Formula for the general number field $K$, which treats all primes in exactly the same way, whether they are discrete or archimedean, and also ramified or not. In another token, I advance a probabilistic interpretation of Weil's positivity criterion, as opposed to the usual geometrical analogies or goals. But in the end, I argue that the new formulation of the Explicit Formula signals a specific link with Quantum Fields, as opposed to the Hilbert-Polya operator idea (which leads rather to Quantum Mechanics).\hfil\par
\vfill
{\parskip = 0 pt
62 rue Albert Joly\par
F-78000 Versailles\par
France\par}

jf.burnol@dial.oleane.com

\eject
\footline{\hfil\the\pageno\hfil}
{\bf TABLE OF CONTENTS}\par
\par
INTRODUCTION\par
FURTHER NOTATIONS AND REVIEW OF WELL-KNOWN RESULTS\par
THE EXPLICIT FORMULA\par
FOURIER TRANSFORM OF $-\log (|X|)$\par
THE CONDUCTOR OPERATOR\par
A PROPAGATOR\par
CONCLUSION\par
REFERENCES\par
\vfill \eject

{\bf INTRODUCTION}\par

The relation between prime powers and the zeros of the zeta function is a striking result of Riemann's paper [1], which can be generalized to $L$-functions and is expressed in various ``Explicit Formulae''. Weil puts forward in his paper from 1952 [2] an identity $$Z(K, \chi )(f) = W(K, \chi )(f)$$ which encompasses most abelian ``Explicit Formulae''.

Here $K$ is a number field, $\chi$ a unitary Dirichlet-Hecke character and $f$ a function on the positive real line, which we will assume to be smooth and compactly supported (although Weil considered functions of a much wider class).

The left-hand-side $Z(K, \chi )(f)$ is obtained by adding the values of the Mellin transform
$$\widehat{f}(s) = \int_{0}^{\infty} f(t){t}^s \, { dt \over t}$$
of $f$ at the poles (counted with positive multiplicity) and the zeros (counted with negative multiplicities) of the $L$-function $L(K, \chi, s)$.

The right-hand-side $W(K, \chi)$ has an additive decomposition  over the prime spots of $K$. Each local term $W_{\nu}(K, \chi)$ is best seen as the push-forward of a distribution $W_{\nu}(K_{\nu}, {\chi}_{\nu})$ from the multiplicative group of the $\nu -$adic field $K_\nu$ to $(0, \infty )$ (under $u_\nu \mapsto t = |u_\nu |_\nu)$. Indeed, all local terms, archimedean as non-archimedean, take then an (almost) identical functional form, and $W(K, \chi )$ is thus best seen as a distribution on the (classes of) ideles of $K$.

All these ``Explicit Formulae'' associated to the various characters $\chi$ are combined by Weil in one unique identity of distributions on the idele classes of $K$. And the Riemann Hypothesis for all the (abelian) $L$-series of the number field $K$ is then shown to be equivalent to a certain ``positivity'' criterion on a related distribution.

We pause here to point out a few minor differences between our conventions and those of [2]: our test-function is on the positive reals whereas Weil's is on the additive reals 	(they are related through a logarithmic transformation); there is a change $s \mapsto {1 \over 2}-s$ from Weil's Mellin transform; the discriminant of the number field appears in [2] as a multiple of the Dirac at $1$, we incorporate it to the local terms; we keep the poles and the zeros together on the same side of the Explicit Formulae.

I said that the local terms took an (almost) identical form, but there are still some differences between the results at archimedean and non-archimedean places. Also, the computation of the local term at an archimedean place is quite involved, and this left its mark in many subsequent expositions.

Haran [3] gave in the case of the Riemann zeta function exactly identical formulations to all local terms. His treatment of the local term at the real prime is more straightforward than Weil's. It is still quite distinct from the computations at the finite places, so that the reason for the coincidence of the functional forms at all places remained mysterious.

For the general number field and character, there arises the problem of ramified primes (primes dividing the discriminant of $K$ or the conductor of $\chi$). At a place dividing the conductor of $\chi$ one needs a specific computation which, perhaps because it is elementary, was in fact not spelled out in [2], only its final result being given. The same dichotomy arises in Weil's later paper [4] which extends the method to the Artin $L$-functions.

In this paper, I will show how to extend Haran's result to the general Dirichlet-Hecke $L$-series (in a reformulated but directly equivalent form). The computations will be exactly the same for all places, whether they are archimedean or not, and also ramified or not. In a later section, and with the positivity in mind, I start exploring some aspects suggested by this reformulation of the Explicit Formula. I point out the relevance of a ``Balayage'' introduced by Zabrodin in a paper on the $p$-adic string [5]. Finally, in the conclusion I hint very briefly on why establishing a deeper link with Quantum Fields appears highly desirable.

\vfill\eject
{\bf FURTHER NOTATIONS AND REVIEW OF WELL-KNOWN RESULTS}

Let $K$ be a number field, $\AA$ its adele ring, $\AA^\times$ its multiplicative idele group, ${\cal C} = {\AA^\times}/{K^\times}$ its multiplicative group of idele classes. Let $F$ be a test-function on $\cal C$, i{.}e{.} a finite linear combination of functions $f(|u|){\chi(u)}^{-1}$ where $|u|$ is the module, $f$ is a smooth function on $(0, \infty)$ with compact support, and $\chi$ is a unitary character on $\cal C$ (the justification for writing ${\chi(u)}^{-1}$ instead of $\chi(u)$ is in our way of defining the Mellin transform).

The Mellin transform of such an $F$ is a function of quasi-characters (continuous homomorphisms of $\cal C$ to $\CC^\times$) defined as
$$\widehat{F}(\chi)=\int_{\cal C}F(u)\cdot\chi(u)\, d^*u$$
The Haar measure $d^*u$ on $\cal C$ in this integral is normalized so that under push-forward to $(0, \infty)$ it is sent to $dt\over t$ (the fibers are compact).

We are now going to recall and use some of the well-known results of Tate's Thesis ([6]). All the Hecke $L$-functions $L(\chi, s)$ can be combined together as one meromorphic function $L(\chi')$ (which can be normalized in various ways) on the space of quasi-characters $\chi'$, according to the formula $L(\chi\cdot\omega_s) = L(\chi, s)$ ($\omega_s(u)=|u|^s$). The multiplicity of $L$ at the quasi-character $\chi'$ will be denoted by $m(\chi')$ (poles counted with positive multiplicities and zeros with negative multiplicities). Let:
$$Z(F)=\sum_{L(\chi')\ =\ 0\ or\ \infty}m(\chi')\cdot\widehat{F}(\chi')$$
This is the left hand side of the Explicit Formula.

Assuming from now on that $F(u) = f(|u|)\chi(u)^{-1}$, for a given fixed unitary character $\chi$, we see that $\widehat{F}(\chi\cdot\omega_s)=\widehat{f}(s)$ while $\widehat{F}(\chi') = 0$ for all other quasi-characters $\chi'$. Writing $Z(\chi, f)$ instead of $Z(F)$ we thus get its value as the sum of $\widehat{f}(s)$ over the zeros and poles of the $L$-function $L(\chi, s)$, as was considered in the Introduction.

It is convenient to consider the $L$-function $L(\chi\cdot\omega_s)$ not as valued in the complex numbers, but as valued in the distributions on the adeles \AA, in the following manner:
$$L(\chi\cdot\omega_s)(\varphi) = \int_{\AA^\times}\varphi(x)\cdot\chi(x)\cdot|x|^s\, d^*x$$
Here the test-function $\varphi$ on the adeles is a finite linear combination of infinite products $\prod_\nu \varphi_\nu$, where $\varphi_\nu$ depends only on the $\nu -$adic component, is for almost all $\nu$'s the characteristic function of the sub-ring of integers in $K_\nu$, is for all non-archimedean places locally constant with compact support (Bruhat function), and is for archimedean places a Schwartz function.

Note that although this defines a distribution on \AA\ the integral is over $\AA^\times$, not over \AA. The Haar measure $d^*x$ is the product of the local multiplicative Haar measures and it is not the same thing as $d^*u$ above which lives on $\cal C$ and was normalized ``globally''. Among other possibilities we normalize the multiplicative Haar measure at an archimedean place so that push-forward to $(0, \infty)$ sends it to $dt \over t$, while at a discrete place we choose it to assign a volume of $1$ to the compact open subgroup of units.

In fact this defines directly $L(\chi\cdot\omega_s)(\varphi)$ only for $Re (s) > 1$ (as an absolutely convergent product). Tate obtained the analytic continuation to all $s$, with the following formula:
$$L(\chi\cdot\omega_s)(\varphi) = \left\{\left\{{{-\kappa\cdot\varphi(0)}\over{s}} + 							{{-\kappa\cdot\widetilde\varphi(0)}\over {1-s}}\right\}\right\}$$
$$+\int_{\AA^{\times}, |x|\geq 1}\{\varphi(x)\cdot\chi(x)\cdot|x|^s + 				\widetilde{\varphi}(x)\cdot{\chi}^{-1}(x)\cdot|x|^{1-s}\}\, d^*x$$
Here $\kappa$ is a constant depending upon the normalizations above and not important to us, $\widetilde{\varphi}(x)$ is the Fourier transform of $\varphi(x)$ and the pole terms within braces appear only when $\chi$ is a principal character ($\chi = \omega_s$ for some purely imaginary $s$).

It follows (``Functional Equation'') that the Fourier transform of the distribution $L(\chi\cdot\omega_s)$ is the distribution $L(\chi^{-1}\cdot\omega_{1-s})$ (note that $\chi(-1) = 1$). The function $L(\chi, s)$ is obtained by evaluating the distribution $L(\chi\cdot\omega_s)$ on a specific choice of test-function, depending on $\chi$ but not on $s$. To discuss this matter further, we need to digress a little through the local situation.

On the local field $K_\nu$, the Fourier transform will be defined using the additive character $\lambda(x) = \exp(2\pi i\{Tr(x)\})$, where $Tr$ is the trace down to $\QQ_p$ (or $\RR$ for an archimedean place, the formula then defines $\lambda(-x)$). There is then a unique choice of additive Haar measure $da$ on $K_\nu$ which is self-dual for $\lambda$. For the global additive Haar measure built from these local choices one finds that the additive group of principal adeles $K \subset \AA$ is its own orthogonal. The (local) module of $x \in K_\nu^\times$ satisfies $d(xa) = |x|_\nu\cdot da$. The multiplicative Haar measure $d^*x$ on $K_\nu^\times$ is a constant multiple of $dx \over |x|_\nu$.

Let $\chi_\nu :K_\nu^\times \rightarrow U(1)$ be a local multiplicative unitary character. For $Re (s) > 0$, $\chi_\nu(x)\cdot|x|_\nu^s$ is locally integrable against $d^*x$ hence defines a distribution on $K_\nu : \varphi_\nu \mapsto L_\nu (\chi_\nu\cdot\omega_s)(\varphi_\nu) = \int_{K_\nu^\times} \varphi_\nu(x) \chi_\nu(x)\cdot|x|_\nu^s\, d^*x$  ($\omega_s$ is here the local principal character $x\mapsto |x|_\nu^s$).

Tate [6] (see also [Weil 7], [Gel'fand 8]) showed that the local distribution $L_\nu (\chi_\nu\cdot\omega_s)$ has a meromorphic continuation to all $s$, and that its Fourier transform is a constant multiple of $L_\nu (\chi^{-1}_\nu\cdot\omega_{1-s})$. The proportionality factor is usually called the $\nu -$adic Tate-Gel'fand-Graev Gamma function and denoted $\Gamma_\nu(\chi_\nu, s)$. It is analytic and non-vanishing in the strip $0 < Re(s) < 1$.

It is tempting to say that the Functional Equation means that the ``adelic'' Gamma Function is identically 1. This is a misleading statement, as we will see now, because the Explicit Formula is basically its logarithmic derivative, and this is not zero!

{\bf THE EXPLICIT FORMULA}

Going back to the global situation, we know from Tate's Thesis that with a suitable test-function $\varphi = \prod_\nu \varphi_\nu$ the local factors of the Hecke $L$-function are obtained as
$$L_\nu(\chi_\nu, s) = L_\nu(\chi_\nu\cdot\omega_s)(\widetilde{\varphi_\nu})$$ (they have neither poles nor zeros in $Re(s) > 0$) while globally
$$L(\chi, s) = L(\chi\cdot\omega_s)(\widetilde\varphi)$$
and is thus obtained for $Re(s) > 1$ as the absolutely convergent product $\prod_\nu L_\nu(\chi_\nu, s)$. For $Re(s) < 0$ one can then represent $L(\chi, s)$ as the absolutely convergent product $\prod_\nu M_\nu(\chi^{-1}_\nu, 1-s)$, where  $M_\nu(\chi^{-1}_\nu, 1-s) = L_\nu(\chi^{-1}_\nu\cdot\omega_{1-s})(\varphi_\nu)$ (it has neither pole nor zero in $Re(s) < 1$). The Tate-Gel'fand-Graev Gamma function is
$$\Gamma_\nu(\chi_\nu, s) = {L_\nu(\chi_\nu, s) \over M_\nu(\chi^{-1}_\nu, 1-s)}$$
and it is without zero or pole in the critical strip $0 < Re(s) < 1$.

The computation of $Z(\chi, f)$ proceeds in the usual way of the calculus of residues as: $$Z(\chi, f) = {1 \over 2\pi i}\int_R \widehat{f}(s)\cdot\ -d\log\ L(\chi,s)$$ where R is the rectangle with corners $2-i\infty, 2 + i\infty, -1+i\infty, -1-i\infty$. For the necessary analytical estimates to discard the horizontal contributions, see [2].

On the vertical line $Re(s) = 2$ one uses the product $\prod_\nu L_\nu (\chi_\nu, s)$ to express the integral as a sum of local terms. As the local $L$-functions have no poles or zeros in the half-plane $Re(s) > 0$, one can shift each local term to the line $Re(s) = c$, where $0 < c < 1$. (That $f$ has compact support turns out to imply that there are only finitely many non-vanishing terms, but this is not necessary to this argument).

On the vertical line $Re(s) = -1$ one uses the product $\prod_\nu M_\nu (\chi_\nu^{-1}, 1-s)$ to express the integral as a sum of local terms. As the local integrands have no poles or zeros in the half-plane $Re(s) < 1$, one can shift each local term to the same line $Re(s) = c$.

Recombining the integrands from the left and the right, we end up with:
$$Z(\chi, f) = \sum_\nu W_\nu(\chi, f)$$
$$\hbox{where }W_\nu(\chi, f) = {1 \over {2\pi i}}\int_{Re(s)=c} \widehat{f}(s)\cdot\Lambda_\nu(\chi, s)\, ds$$
$$\Lambda_\nu(\chi, s) = - \hbox{logarithmic derivative of }\Gamma_\nu(\chi, s)$$

This already differs from Weil's method in [2] in the way one uses the functional equation but has the advantage of treating all places, ramified or not, alike. The next classical step from then on would be to obtain the inverse Mellin transform $w_\nu(\chi, t)$ of $\Lambda_\nu(\chi, s)$ (in the critical strip), which is a distribution on the positive half-line such that
$$W_\nu(\chi, f) = \int_0 ^\infty f({1 \over t})w_\nu(\chi, t)\, {dt \over t}$$
and then (somewhat surprisingly in the real case and quite surprisingly in the complex case) to realize that it comes from $K_\nu^\times$ through $u \mapsto t = |u|_\nu$\ .

Let us now work locally and drop all subscripts $\nu$ to lighten the notation. Let $\cal F$ be the operation of Fourier transform on the $\nu$-adic completion of $K$. The identities to follow are identities of locally integrable functions, taken in the distributional sense, and the parameter $s$ is in the critical strip $0<Re(s)<1$.
$${\cal F}(\chi(x)|x|^{s-1}) = \Gamma(\chi, s)\chi^{-1}(y)|y|^{-s}$$
$${\cal F}\bigl(\log|x|\,\chi(x)|x|^{s-1}\bigr)=\Gamma^\prime(\chi,s)\chi^{-1}(y)|y|^{-s} + \Gamma(\chi, s)\,(-\log|y|)\chi^{-1}(y)|y|^{-s}$$
$$\chi^{-1}(y)\Lambda(\chi,s)|y|^{-s}=-\log|y|\,\chi^{-1}(y)|y|^{-s} + {\cal F}\bigl(-\log|x|\cdot\chi(x){|x|^{s-1}\over\Gamma(\chi, s)}\bigr)$$
We think of these identities as being applied to some given test-function and we now want to apply a further operation of integration, this time against ${1 \over {2\pi i}}\int_{Re(s)=c} \widehat{f}(s)\cdot\ \, ds$, for $0<c<1$. One needs a lemma according to which (pointwise)
$${1 \over {2\pi i}}\int_{Re(s)=c} \widehat{f}(s)\cdot\chi(x){|x|^{s-1}\over\Gamma(\chi, s)}\, ds
= {\cal F}^{-1}(\chi^{-1}(y)f(|y|))$$
The integral converges absolutely due to the rapid decrease of $\widehat{f}(s)$ whereas $\left|{1\over\Gamma(\chi,s)}\right|= \left|\Gamma(\chi^{-1},1-s)\right|$ ($=1$ on the critical line) is in any case periodic for a finite place, and controlled by the Stirling formula to be $O(|s|^{1/2-c})$ for a real place, $O(|s|^{1-2c})$ for a complex place. So the left-hand-side defines a continuous function of $x$ (for $x\neq0$) and also a tempered distribution as it is $O(|x|^{c-1})$. Applying both sides to a test function we see  (using Mellin inversion) that they coincide as distributions hence pointwise (for $x\neq0$).

As everything converges absolutely we are allowed to intervert integrals and end up with
$$W_\nu(\chi, f; y) = -\log(|y|_\nu)\cdot F_\nu(y) + {\cal F}\bigl(-\log(|x|_\nu)\cdot  {\cal F}^{-1}(F_\nu(y))\bigr)$$
where we have put back the $\nu$ subscripts, defined $F_\nu(y) = f(|y|_\nu)\cdot\chi_\nu^{-1}(y)$ to be the local component of $F$ and (for $y \in K_\nu^\times$):
$$W_\nu(\chi, f; y) = {1 \over {2\pi i}}\int_{Re(s)=c} \widehat{f}(s)\cdot\Lambda_\nu(\chi, s)\cdot|y|_\nu^{-s}\cdot\chi_\nu^{-1}(y)\, ds$$
$$W_\nu(\chi, f; y) = -\log(|y|_\nu)\cdot F_\nu(y) + {\cal F}\bigl(-\log(|x|_\nu)\bigr) * F_\nu(y)$$
where $*$ is the symbol of additive convolution on $K_\nu$ . Plugging in $y = 1$ gives:
$$W_\nu(F) = W_\nu(\chi, f) = (G_\nu * F_\nu)(1)$$
where $G_\nu$ is the Fourier Transform of $-\log(|x|_\nu)$ and $*$ is an additive convolution. For $K = \QQ_p$, or $\RR$, and with $\chi$ the trivial character, this is directly equivalent to a result of Haran [3], and our goal was to extend his formula to the general case. In doing so, it turned out to be possible to give a uniform derivation for all places of the arbitrary number field $K$.

Here is the Explicit Formula:
$$\sum_{L(\chi')\ =\ 0\ or\ \infty}m(\chi')\cdot\widehat{F}(\chi') = \sum_\nu(G_\nu * F_\nu)(1)$$
$$\hbox{for an arbitrary test-function\ }F\hbox{\ on the idele classes,}$$
$$\hbox{with\ }F_\nu \hbox{\ the restriction of\ } F \hbox{\ to\ }K_\nu^\times$$
$$\hbox{and\ }G_\nu = {\cal F}(-\log(|x|_\nu))$$

{\bf FOURIER TRANSFORM OF $-\log(|X|)$}

To put the result in Weil's form we need to obtain the Fourier transform of $-\log(|x|)$. This is a classical calculus exercise in the archimedean case (see, for example, Lighthill [9]). Over $\QQ_p$ it has been obtained by Vladimirov [10]. It is also implicitely contained in Haran's work. As previously for the Explicit Formula we derive it in a simple manner from the properties of homogeneous distributions considered in [Gel'fand 8], [Weil 7].

Again we work locally, and will drop most $\nu$ subscripts to ease up on the notations. Let $\Delta_s$ be the homogeneous distribution, given for $Re(s) > 0$, by $\Delta_s(\varphi) = \int_{K_\nu}\varphi(x)\cdot|x|^{s-1}\, dx$. This time we need an explicit expression for its analytic continuation around $s = 0$. Let us (temporarily) pick a specific function $\omega(x)$ as follows:

If $\nu$ is finite, we take $\omega$ to be the characteristic function of the ring of integers $\cal O$. Let $\pi$ be a uniformizer and $q = {|\pi|}^{-1}$ (which is also the cardinality of the residue field). Let $\delta$ be the differental exponent at $\nu$, which is the largest integer such that $(x \in \pi^{-\delta}{\cal O}) \Rightarrow (Tr(x) \in \ZZ_p)$. Then Vol$({\cal O}) = q^{(-\delta/2)}$, and:
$$\Delta_s(\omega) = q^{-(\delta/2)} {(1 - 1/q) \over (1 - 1/q^s)}$$
So $\Delta_s(\omega)$ has a meromorphic continuation to all $s$, and around $0$:
$$\Delta_s(\omega) = {R\over s} + P(\omega) + \cdots$$
with $R = q^{-(\delta/2)} {(1 - 1/q) \over \log(q)}$ and $P(\omega) = 1/2\cdot\log(q)\cdot R$

If $\nu$ is real, we take $\omega(x) = \exp(- \pi x^2)$ and then:\hfil\break
\centerline{$\Delta_s(\omega) = \pi^{-s/2}\, \Gamma(s/2)$ (here $\Gamma$ is the classical Gamma function)}
\centerline{around $0$: $\Delta_s(\omega) = {R\over s} + P(\omega) + \cdots$}
with $R = 2$ and $P(\omega) = -(\log(\pi) + \gamma_e)$ ($\gamma_e = -\Gamma^\prime(1)$ = Euler's constant)

If $\nu$ is complex we take $\omega(z) = \exp(- \pi z\overline z)$ and recall that $|z|$ means $z\overline z$, and that the self-dual Haar measure for $\lambda(z) = \exp(-2\pi i\ 2Re(z))$ is twice the Lebesgue measure, so that $\Delta_s(\omega) = 2\pi\ \pi^{-s}\, \Gamma(s)$:
$$\Delta_s(\omega) = {R\over s} + P(\omega) + \cdots$$
with $R = 2\pi$ and $P(\omega) = - 2\pi\cdot(\log(\pi) + \gamma_e)$

The continuation to $Re(s) > -1$ (to all $s$ if $\nu$ is non-archimedean) of $\Delta_s$ is obtained (for example) by:
$$\Delta_s(\varphi) = \Delta_s(\omega)\cdot\varphi(0) + \int_{K_\nu}(\varphi(x)-\varphi(0)\omega(x))\, |x|^s\, {dx \over |x|}$$
One sees that $\Delta_s$ has a simple pole at $0$ ($\delta$ = Dirac at $0$):
$$\Delta_s = {{R\cdot\delta}\over s} + P + \cdots$$
with the ``finite part'' $P(\varphi) = P(\omega)\varphi(0) + \int_{K_\nu}(\varphi(x)-\varphi(0)\omega(x))\, {dx \over |x|}$ is expressed in terms of $\omega$ but is independant of it.

This shows that we could have taken any other test-function $\omega$ such that $\omega(0) = 1$ and, adopting the notation $P_\omega$ for the distribution $\int_{K_\nu}(\varphi(x)-\varphi(0)\omega(x))\, {dx \over |x|}$, that $P = P_\omega + P(\omega)\cdot\delta$.

We know that ${\cal F}(\Delta_s) = \Gamma(s)\cdot\Delta_{1-s}$ with $\Gamma(s)$ the Gel'fand-Graev Gamma function for the trivial character. Expanding this in powers of $\varepsilon$ with $s = 1 - \varepsilon$, we get, on the left hand side:
$${\cal F}(\Delta_{1-\varepsilon}) =  \delta + \varepsilon\cdot{\cal F}(-\log|x|) + \cdots$$
On the right hand side, we have, introducing constants $\gamma$ and $\tau$:
$$\Gamma(1 - \varepsilon) = \gamma\cdot(\varepsilon + \tau\cdot\varepsilon^2 + \cdots)$$
$$\Delta_\varepsilon = {{R\cdot\delta}\over\varepsilon} + P + \cdots$$
Hence $\gamma\cdot R = 1$ and the formula for $G = {\cal F}(-\log|x|)$ is:
$$G = \gamma\cdot P + \tau\cdot\delta$$
$$G = \gamma\cdot P_\omega + (\gamma\cdot P(\omega) + \tau)\cdot\delta$$
$$G = \gamma\cdot P_\omega + G(\omega)\cdot\delta$$
$G(\varphi) = \gamma\cdot P(\varphi) + \tau\cdot \varphi(0)$ is also directly given by:
$$G(\varphi) =  \int_{K_\nu}(-\ln|y|)\, \widetilde{\varphi}(y)\, dy = - \left.{\partial \over {\partial s}}\right|_{s = 1}\ \Delta_s(\widetilde{\varphi})$$
from which follows for $\varphi_t(x) = \varphi(tx)$:
$$G(\varphi_t) = G(\varphi) - \log(|t|)\cdot\varphi(0)$$
For all places, the constant $\gamma = 1/R = -\Gamma^\prime(1)$ turns out to be such that the multiplicative Haar measure $d^\times u = \gamma\cdot {du \over |u|}$ on $K_\nu^\times$ has the property  $\int_{1\leq|x|<|X|}\, d^\times u = \log(|X|)$ (with $X \in K_\nu, |X| \geq 1)$. For a discrete $\nu$ it assigns a value of $\log(q)$ to the volume of the units.

Let us suppose now that the place $\nu$ is non-archimedean. For $\omega$ the characteristic function of the integers, its Fourier transform $\widetilde \omega$ is $\widetilde\omega(y) = q^{-\delta/2}\cdot\omega(\pi^\delta\cdot y)$ and $\Delta_s(\widetilde\omega) = q^{\delta\cdot(s-1)}\cdot{(1 - 1/q) \over (1 - 1/q^s)}$ so that $G(\omega) = -\left.{\partial \over {\partial s}}\right|_{s = 1}\ \Delta_s(\widetilde{\omega}) = {\log(q)\over (q-1)} - \delta\cdot\log(q)$. So:
$$G(\varphi)= \gamma\cdot(P_\omega(\varphi) + q^{-(\delta/2)-1}\cdot\varphi(0)) - \delta\cdot\log(q)\cdot\varphi(0)$$

Let us apply this to the characteristic function $\psi$ of the ball of radius $1/q$ centered at $1$: we get $G(\psi) = \gamma\cdot(q^{-(\delta/2)-1})$. So it is advantageous to replace $\omega$ by $\omega_1 = \omega - \psi$, and then:
$$G(\varphi) = \gamma\cdot P_{\omega_1}(\varphi) - \delta\cdot\log(q)\cdot\varphi(0)$$
Going back to Weil's local term $W_\nu(F) = (G_\nu * F_\nu)(1)$, this gives:
$$W_\nu(F) =  \int_{K_\nu}(F_\nu(1-x) - F_\nu(1)\omega_1(x))\, \gamma\cdot{dx \over |x|} - \delta\cdot\log(q)\cdot F_\nu(1)$$
and after the change of variable $x \mapsto 1 - 1/u$, we get:
$$W_\nu(F) = \int_{|u| = 1}(F_\nu({1\over u}) - F_\nu(1))\, {d^\times u \over |1 - u|} + \int_{|u| \neq 1}F_\nu({1\over u})\, {d^\times u \over |1 - u|}-\ \delta\cdot\log(q)\cdot F_\nu(1)$$
This is Weil's local term ([2]), with the addition of the local component of the discriminant. The integral over $|u| = 1$ for the test-function $F(u) = f(|u|)\cdot\chi^{-1}(u)$ gives a non-zero result only if $\chi$ is ramified at $\nu$, whereas the other integral gives a non-zero result only when $\chi$ is non-ramified at $\nu$. We will come back to the ramified characters, but first let us examine the local term at an archimedean spot.

First we consider a real spot $\nu$.\hfil\break For $\omega(x) = \exp(- \pi\ x^2)$ we evaluated $\Delta_s(\omega) = \pi^{-s/2}\, \Gamma(s/2)$, $R = 2$ (so $\gamma = 1/2$), and $P(\omega) = -(\log(\pi) + \gamma_e)$. Using $\omega = \widetilde\omega$  (and also $\Gamma^\prime(1/2) = -\ \Gamma(1/2)\cdot(\gamma_e + 2\ \log(2)) )$ we get $G(\omega) = - \left.{\partial \over {\partial s}}\right|_{s = 1}\ \Delta_s(\widetilde{\omega}) = (\log(\pi) + \gamma_e + 2\log(2))/2$. From this we get the constant $\tau = G(\omega) - \gamma\cdot P(\omega) = \log(2\pi) + \gamma_e$.

It is convenient now to choose other functions $\omega$ (always obeying the condition $\omega(0) = 1$). For example, we can take $\omega(x) = 1$ if $|x| < 1$, $\omega(x) = 0$ if $|x| \geq 1$. (it doesn't matter that $\omega$ is not smooth away from $0$). For this $\omega$, $\Delta_s(\omega) = 2/s$, so $P(\omega) = 0$ and $G(\omega) = \tau = \log(2\pi) + \gamma_e$. So we get:
$${\cal F}(-\log|x|)(\varphi) = \int_{|x| \leq 1}(\varphi(x) - \varphi(0))\, {dx \over 2|x|} + \int_{|x| > 1}\varphi(x)\, {dx \over 2|x|}+\ (\log(2\pi) + \gamma_e)\cdot\varphi(0)$$
The local Weil term $W_\RR(F)$ is then obtained as:
$$\int_{1/2}^\infty(F({1\over u}) - F(1))\, {d^\times u \over |1 - u|} + \int_{-\infty}^{1/2}F({1\over u})\, {d^\times u \over |1 - u|} +\ (\log(2\pi) + \gamma_e)\cdot F(1)$$
with the multiplicative Haar measure $d^\times u = {du \over 2|u|}$.

Choosing $\omega(x) = \sin(\pi x)/\pi x$ whose Fourier transform (as a distribution) is the characteristic function of $-1/2 \leq x \leq 1/2$, we get $\Delta_s(\widetilde\omega) = 1/s\cdot(1/2)^{s-1}$ hence $G(\omega) = 1 + \log(2)$ and:
$${\cal F}(-\log|x|)(\varphi) =  \int_\RR \left({\varphi(x) - \varphi(0){\sin(\pi x) \over \pi x}}\right)\, {dx \over 2|x|} + (1 + \log(2))\cdot\varphi(0)$$
This leads to the following expression for Weil's term:
$$W_\RR(F) = \int_{\RR^\times}\left({F({1\over u}) - {{{u\over\pi}\, \sin({\pi\over u})}\over u - 1}\, F(1)}\right)\, {d^\times u\over|1 - u|} + (1 + \log(2))\cdot F(1)$$
As $\exp(G(\omega)) = 2e$, the function $\omega(x) = \sin(2\pi ex)/2\pi ex$ would lead to a similar expression with no Dirac term.

Choosing $\omega(x) = \left({\sin(\pi x)/\pi x}\right)^2$ whose Fourier transform is the triangle function $x \mapsto 1 - |x|$ for $|x| \leq 1$, we get $\Delta_s(\widetilde\omega) = 2/s(s+1)$ hence $G(\omega) = 3/2$ and:
$${\cal F}(-\log|x|)(\varphi) =  \int_\RR \left({\varphi(x) - \varphi(0){\left({\sin(\pi x) \over \pi x}\right)^2}}\right)\, {dx \over 2|x|} + (3/2)\cdot\varphi(0)$$
$$W_\RR(F) = \int_{\RR^\times}\left({F({1\over u}) - {\left({{{u\over\pi}\, \sin({\pi\over u})}\over u - 1}\right)^2}\, F(1)}\right)\, {d^\times u\over|1 - u|} + (3/2)\cdot F(1)$$

Switching to a complex spot $\nu$, we obtained with $\omega(z) = \exp(-\ \pi z\overline z)$: $\Delta_s(\omega) = 2\pi\ \pi^{-s}\, \Gamma(s)$, $R = 2\pi$, $\gamma = 1/2\pi$ and $P(\omega) = -2\pi\cdot(\log(\pi) + \gamma_e)$. The Fourier transform (recall that $dz = 2rdrd\theta$ and $|z| = z\overline z$) of $\omega(z)$ is $2\omega(2z)$, so $\Delta_s(\widetilde \omega) = (4\pi)^{1-s}\ \Gamma(s)$ and we get $G(\omega) = - \left.{\partial \over {\partial s}}\right|_{s = 1}\ \Delta_s(\widetilde{\omega}) =  \log(4\pi) + \gamma_e$. From this we get the constant $\tau = G(\omega) - \gamma\cdot P(\omega) = 2(\log(2\pi) + \gamma_e)$. Choosing $\omega(z)$ to be the characteristic function of the unit disc, we find $\Delta_s(\omega) = 2\pi/s$ hence $P(\omega) = 0$ and so $G(\omega) = \tau$, and the following thus emerges:
$${\cal F}(-\log|x|)(\varphi) = \int_{|x| \leq 1}(\varphi(x) - \varphi(0))\, {dx \over 2\pi|x|} + \int_{|x| > 1}\varphi(x)\, {dx \over 2\pi|x|}+\ 2\cdot(\log(2\pi) + \gamma_e)\cdot\varphi(0)$$
This gives the following expression for the Weil term $W_\CC(F)$ ($d^\times u = drd\theta/\pi r$):
$$\int_{Re(u)\geq 1/2}(F({1\over u}) - F(1))\, {d^\times u \over |1 - u|} + \int_{Re(u)\leq 1/2}F({1\over u})\, {d^\times u \over |1 - u|}+\ 2\cdot(\log(2\pi) + \gamma_e)\cdot F(1)$$

There are of course countless other possible choices for the function $\omega$ and ensuing expressions for $W_\RR(F)$ and $W_\CC(F)$.

{\bf THE CONDUCTOR OPERATOR}

We now return to a discrete place $\nu$ of the number field $K$, with $\nu$-adic completion $K_\nu$. Let $\chi$ be a local, ramified, unitary character (ramified meaning that the restriction of $\chi$ to ${\ |u| = 1\ }$ is non-trivial). Let $f\geq1$ be the conductor exponent of $\chi$, that is the smallest positive integer such that $| 1 - u | \leq q^{-f} \Rightarrow \chi(u) = 1$. For $q=2$ the smallest possible value of $f$ is $2$ whereas for all other $q$'s it is $1$. The Gamma function of $\chi$ was computed by Tate to be simply $q^{(f+\delta)\cdot s}$ up to a certain multiplicative constant, not important here. So minus its logarithmic derivative is the constant $-(f + \delta) \log(q)$.

Going back to the evaluation of $$W_\nu(\chi, f; y) = {1 \over {2\pi i}}\int_{Re(s)=c} \widehat{f}(s)\cdot\Lambda_\nu(\chi, s)\cdot|y|_\nu^{-s}\cdot\chi_\nu^{-1}(y)\, ds$$
this gives $W_\nu(\chi, f; y) = -(f + \delta) \log(q)\cdot f(|y|)\cdot \chi^{-1}(y)$, and for $F(y) = f(|y|)\cdot\chi^{-1}(y)$ our previous computation of $W_\nu(\chi, f; y)$ then gives:
$$-(f + \delta) \log(q) F(y) = (-\log|y|)\cdot F(y) + {\cal F}(-\log|x|) * F(y)$$
Choosing the function $f$ to have its support contained in $1/q < u < q$, and such that $f(1) = 1$, plugging in $y = 1$, and using our evaluation of ${\cal F}(-\log|x|)$, we end up with:$$-(f + \delta) \log(q) =  \int_{|u| = 1}(\chi(u) - 1){{d^\times u}\over{|1-u|}} - \delta\cdot\log(q)$$
$$f \log(q) =\int_{|u| = 1}(1 - \chi(u)){{d^\times u}\over{|1-u|}}$$
It is not difficult to evaluate directly this integral which is a necessary component of the final result of [2] (perhaps for this reason a proof was omitted there).

Let us give an operator interpretation for this ramified case of the local term. Let $\cal L$ be the $L^2$ Hilbert space with respect to the additive Haar measure $dx$. We have various unitary operators on $\cal L$:
\+&the Fourier transform $\cal F$\cr
\+&the Rotations and Dilations: ${\cal R}(u): F(x) \mapsto |u|^{1/2} F(ux)$\cr
\+&the Inversion ${\cal I}: F(x) \mapsto 1/|x| F(1/x)$\cr
\+satisfying the commutation relations:\cr
\+&${\cal R}(u){\cal F} = {\cal F\ R}(1/u)$\cr
\+&${\cal R}(uv) = {\cal R}(u){\cal R}(v)$\cr
\+and\hfil&${\cal R}(u){\cal I} = {\cal I\ R}(1/u)$\cr
They act on the (``cuspidal''?) subspace ${\cal L}_0$ which is defined as the null-space of the averaging operator $\int_{|u| = 1}{\cal R}(u)\, d^*u$ and is spanned by functions $\varphi(x) = g(|x|)\chi(x)$, where $g$ has compact support on $(0, \infty)$ and $\chi$ is a ramified character. Let $A$ be the operator with dense domain of definition defined by $A(\varphi)(x) = \log(|x|)\cdot \varphi(x)$. The spectrum of $A$ has its support consisting of all integer multiples of $\log(q)$, of both signs. We have the following commutation relations:
\+&${\cal R}(u)\ A = A\ {\cal R}(u) + \log(|u|)\ {\cal R}(u)$\cr
\+&${\cal I}\ A = - A\ {\cal I}$\cr
$A$ does not commute with the Fourier transform, so let $B = {\cal F}A{\cal F}^{-1}$ and $H = A + B$. The commutation relations above give:\par
\+&${\cal R}(u)\ B = B\ {\cal R}(u) - \log(|u|)\ {\cal R}(u)$\cr
\+so\hfill&${\cal R}(u)\ H = H\ {\cal R}(u)$\cr

As $\log|x|$ is even $H$ commutes with the Fourier transform. We call $H$ the ``conductor operator''. In fact we have already computed an explicit formula for $H$. Applied to $\varphi(x) = g(|x|)\chi(x)$ it just gives $(f + \delta) \log(q) \varphi(x)$, so we have the complete spectral decomposition of $H$:

$H$ is a self-adjoint, positive operator, whose spectrum has its support consisting of the integer multiples $(f + \delta) \log(q)$ of $\log(q)$, with $f \geq 1$ for $q>2$, $f\geq2$ for $q=2$, and the eigenspace corresponding to this eigenvalue is the closure of functions $g(|x|)\chi(x)$, where $g$ has compact support and $\chi$ is a character with conductor exponent $f$.

$H$ commutes with the Fourier transform, with the Inversion, and with the operators ${\cal R}(u)$. Also it appears now that $H$ commutes with $A$, hence with $B$, and also that $A$ and $B = H - A$ commute (on a dense subdomain of the cuspidal space ${\cal L}_0$).

{\bf A PROPAGATOR}

Let us now start exploring what the new functional form given to the Explicit Formula suggests about the sought-for positivity.
$$\sum_{L(\chi')\ =\ 0\ or\ \infty}m(\chi')\cdot\widehat{F}(\chi') = \sum_\nu(G_\nu * F_\nu)(1)$$
First of all, the famous aphorism of Einstein about his equation of general relativity comes to mind: the Explicit Formula is like a castle with two sides; the left-hand-side $Z(F)$ is in marble while the right-hand-side $W(F)$ (rather its decomposition in local terms) is in plain wood. Even in the new ``$-\log(|x|)$'' formulation there are arbitrary choices. The following alternative expression is equally valid:
$$\sum_{L(\chi')\ =\ 0\ or\ \infty}m(\chi')\cdot\widehat{F}(\chi') = \sum_\nu(G_\nu * F_\nu^{(k)})(k)$$
Here $k$ is an arbitrary non-vanishing element of the number field $K$ and $$F_\nu^{(k)}(x_\nu) = F((k, k, k,\dots, x_\nu, k, k,\dots))$$ while $F_\nu(x_\nu)$ was defined as $F((1, 1, 1,\dots, x_\nu, 1, 1,\dots))$. The two decompositions are related by:
$$(G_\nu * F_\nu^{(k)})(k) = (G_\nu * F_\nu)(1) + \log(|k|_\nu)\cdot F(1)$$
So in a way the multiplicative group $K^\times$ is a symmetry group of the Explicit Formula.

But of course the analogy with Einstein's remark really comes into being when one recalls that he spent decades trying to give a geometric interpretation to the wood side of his equation, and that arithmeticians too seem to have been pursuing during these last decades a similar goal for $W(F)$.

I will make a few remarks here. To express his positivity criterion Weil uses conventions slightly distinct from ours. He moves the local terms to be together with the poles, and makes a shift of $1/2$ in the Mellin transform. In this way he gets a distribution $C$ and translates the Riemann Hypothesis into a positivity criterion:
$$C(F*F^\tau) \geq 0$$
\centerline{for an arbitrary test-function $F$ on $\cal C$}
\centerline{$F^\tau(x) = \overline{F({1 \over x})}$ and $*$ is a multiplicative convolution}

In the function field case $C(F*F^\tau)$ can be given a geometric interpretation as an intersection number of cycles on an algebraic surface, and the positivity follows from the Hodge Index Theorem (one reference on this story is Haran [11]).

But another interpretation is possible that does not seem to have been pushed forward so far. To prove that a number is non-negative it is enough to exhibit it as the variance of a random variable. In our case this means that there should be a generalized, stationary, zero mean, stochastic process with $\cal C$ as ``time'' whose covariance would be $C$. That is we have a probability measure $\mu$ on the distributions on the classes of ideles and the identity:
$$\int_{\hbox{Distributions}}X_F(d)\overline{X_G(d)}\, d\mu(d) = C(F*G^\tau)$$
where $F$ and $G$ are two arbitrary test-functions on $\cal C$ and $X_F$ and $X_G$ are the two associated ``coordinates'' on the space of distributions $d \mapsto X_F(d) = d(F)$, $d \mapsto X_G(d) = d(G)$. So $C(F*G^\tau)$ is the covariance of the random variables defined by $F$ and $G$. For $F = G$ it is obviously non-negative!

If such a probability measure $\mu$ could be constructed, corresponding to an ``arithmetic stochastic process'', then the Riemann Hypothesis would follow of course.

But such a generalized process could not be obtained and I now believe that this ``temporal'' interpretation of the idele classes and the accompanying purely probabilistic formulation of the Riemann Hypothesis are slightly misguided. We should keep the poles and the zeros together and not make the shift by $1/2$.

There is an essential concept from Quantum Physics which I believe will play an important role in this problem. This is the notion of correlation functions of a quantum field. Well-known advances in the field of knot invariants and many other geometric arenas have shown that these correlation functions can be used to represent intersection products or linking numbers and many other geometrical things. On the other hand they have of course a probabilistic interpretation (being complex amplitudes) and rigorous mathematical developments of constructive quantum field theory have as their goal the construction of probability measures on spaces of distributions, as considered above.

For these reasons, and other more precise reflections on the subject, moving towards quantum fields is an urgent goal. In this context it is reassuring that the ``$-\log(|x|)$'' formulation of the Explicit Formula enables a few additional comments.

Indeed as is well-known ``$-\log(|x|)$'' is the propagator of the free Boson field in 2 dimensions (here we look at the complex place and $|x| = x\overline x$ so that we have both holomorphic and anti-holomorphic sectors). But let us rather consider the non-archimedean completion $K_\nu$ of the number field $K$. Here too ``$-\log(|x|)$'' can be seen as a propagator associated to an action. The main difference with the complex case being that this action turns out to be non-local in position space. But it has a simple Dirac-like (or rather ``square-root of Laplacian''-like) expression in momentum space:
$$S(\varphi) = {q^{-\delta /2}\cdot(1 - {1 \over q}) \over \log(q)}\cdot{1 \over 2}\int_{K_\nu}|k|\varphi(k)\varphi(-k)\, dk$$

\centerline{(for simplicity $\varphi$ is assumed real-valued in position space)}
\centerline{(the numerical factor is the constant $R = 1/\gamma$)}
The usual recipe would give $\gamma/|k|$ as the propagator (in momentum space) but there is no distribution with this homogeneity on $K_\nu$ apart from a Dirac at the origin. An infra-red cut-off is needed and we end up with one of the distributions $\gamma\cdot P_\omega$ considered before, up to a Dirac at the origin. The propagator in position space will be its Fourier transform, which we know to be $-\log(|x|)$, up to an overall additive constant. This is the freedom we observed before in the Explicit Formula.

In position space the action is non-local and skipping a few calculational steps, we get: (the numerical factor in front of the integral sign is $-1/4\cdot R\cdot\Gamma(2)$)
$$S(\varphi) = {{q\cdot(q-1)\cdot q^\delta}\over{4\cdot(q+1)\cdot\log(q)}}\cdot\int_{K_\nu\times K_\nu}{{(\varphi(x)-\varphi(y))^2}\over{|x-y|^2}}\, dx\,dy$$
Whether looking at this in position or in momentum space we are faced with rather non-inspiring numerical factors in front of these integrals. But they have been considered before by physicists especially in the context of $p$-adic strings. In particular, Zabrodin [5] has obtained an inspiring result (at least for $K_\nu = \QQ_p$):  this action with exactly all its numerical factors is obtained as the swept-out to the boundary (Poincar\'e Balayage) action from a simple Gaussian model on the Bruhat-Tits tree with boundary $\PP^1(\QQ_p)$.

Possible links between the Riemann Hypothesis and $SL(2)$ or $PGL(2)$ have been advocated, I think that the result of Zabrodin could play an important role in such a context. But the positivity has remained elusive so far.

{\bf CONCLUSION}

This paper finds its roots in the conviction that the Riemann Hypothesis has a lot to do with (suitably envisioned) Quantum Fields. The belief in a possible link between the Riemann Hypothesis and Quantum Mechanics seems to be widespread and is a modern formulation of the Hilbert-Polya operator approach. I believe that techniques and philosophy more organic to Quantum Fields will be most relevant. As this point of view has not so far led to success, I will conclude here with extracts from William Blake's ``The Marriage of Heaven and Hell'' which, to my mind, remarkably illustrate some salient points :\par
{\parskip = 0 pt
\font\smallerslanted = cmsl8
\smallerslanted
``[\dots]\par
But first the notion that man has a body distinct from his soul is to be expunged; this I shall do, by printing in the infernal method, by corrosives, which in Hell are salutary and medicinal, melting apparent surfaces away, and displaying the infinite which was hid.\par
If the doors of perception were cleansed every thing would appear to man as it is, infinite.\par
For man has closed himself up, till he sees all things thro' narow chinks of his cavern.\par
[\dots]\par
By degrees we beheld the infinite Abyss, fiery as the smoke of a burning city; beneath us at an immense distance, was the sun, black but shining; round it were fiery tracks on which revolv'd vast spiders, crawling after their prey; which flew, or rather swum, in the infinite deep, in the most terrific shapes of animals sprung from corruption; \& the air was full of them, \& seem'd composed of them: these are Devils, and are called Powers of the air. I now asked my companion which was my eternal lot? he said, `between the black \& white spiders.'\par
[\dots] \par
Opposition is true Friendship.\par
[\dots]''}
\vfill \eject
{
{\bf REFERENCES}\par
\baselineskip = 12 pt
\parskip = 4 pt
\font\smallRoman = cmr8
\smallRoman
\font\smallBold = cmbx8
\font\smallSlanted = cmsl8
{\smallBold [1] B. Riemann},{\smallSlanted ``\"Uber die Anzahl der Primzahlen unter einer gegebenen Gr\"osse''}, Monatsber. Akadem. Berlin, 671-680, (1859).

{\smallBold [2] A. Weil},{\smallSlanted ``Sur les ``formules explicites'' de la th\'eorie des nombres premiers''}, Comm. Lund (vol d\'edi\'e \`a Marcel Riesz), (1952).

{\smallBold [3] S. Haran},{\smallSlanted ``Riesz potentials and explicit sums in arithmetic''}, Invent. Math.  101, 697-703 (1990).

{\smallBold [4] A. Weil},{\smallSlanted ``Sur les formules explicites de la th\'eorie des nombres''}, Izv. Mat. Nauk. (Ser. Mat.) 36, 3-18, (1972).

{\smallBold [5] A. V. Zabrodin},{\smallSlanted ``Non-Archimedean Strings and Bruhat-Tits Trees''}, Commun. Math. Phys. 123, 463-483 (1989).

{\smallBold [6] J. Tate}, Thesis, Princeton 1950, reprinted in Algebraic Number Theory, ed. J.W.S. Cassels and A. Fr\"ohlich, Academic Press, (1967).

{\smallBold [7] A. Weil},{\smallSlanted ``Fonctions z\^etas et distributions''}, S\'eminaire Bourbaki n${}^{\smallRoman o}$ 312, (1966).

{\smallBold [8] I. M. Gel'fand, M. I. Graev, I. I. Piateskii-Shapiro},{\smallSlanted ``Representation Theory and automorphic functions''}, Philadelphia, Saunders (1969).

{\smallBold [9] M. J. Lighthill},{\smallSlanted ``Introduction To Fourier Analysis And Generalised Functions''}, Cambridge University Press (1959).

{\smallBold [10] V. S. Vladimirov},{\smallSlanted ``Generalized functions over the fields of p-adic numbers''}, Russian Math. Surveys 43:5, 19-64 (1988).

{\smallBold [11] S. Haran},{\smallSlanted ``Index Theory, Potential Theory, and the Riemann Hypothesis''}, in {\smallSlanted ``L-functions and Arithmetic, Durham Symposium''}, LMS Lecture Note Series 153, eds. Coates \& Taylor, Cambridge University Press (1989).

{\smallBold [12] W. Blake},{\smallSlanted ``The Marriage of Heaven and Hell''}, engraved and coloured by the author, London (1793).

\vfill
\centerline{Jean-Fran\c{c}ois Burnol}
\centerline{62 rue Albert Joly}
\centerline{F-78000 Versailles}
\centerline{France}
\centerline{jf.burnol@dial.oleane.com}
\centerline{September 1998}
\centerline{revised November 1998}
}
\eject
\bye